
\input amstex
\documentstyle{amsppt}
\pagewidth{125mm}
\pageheight{215mm}

\topmatter
\title
Periodic orbits in magnetic fields and Ricci curvature of
Lagrangian systems
\endtitle
\author
Abbas BAHRI and Iskander A. TAIMANOV
\endauthor
\address
Department of Mathematics, Rutgers University, New Brunswick,
New Jersey 08903, USA
\endaddress
\address
Institute of Mathematics, 630090 Novosibirsk, Russia
\endaddress
\email
taimanov\@math.nsk.su
\endemail
\endtopmatter
\leftheadtext{}
\rightheadtext{}
\document
\NoBlackBoxes

\head
\S 1. Introduction
\endhead

Let $M^n$ be a differentiable $n$-dimensional closed manifold endowed
with a Riemannian metric $g_{ij}$ and with an exact 2-form $F$.

We will consider the problem of existence of closed extremals of the
functional
$$
S(\gamma)=\int_{\gamma}(\sqrt{g_{ij}\dot x^i\dot x^j} + A_i\dot x^i)dt,
$$
where $dA=F$, on the space of closed curves on the manifold $M^n$, and
where $A$
is a 1-form (i.e., $F$ is
an exact 2-form).

This functional is a natural generalization of the usual functional of
length, and its closed extremals correspond to periodic trajectories of the
motion of particle on the Riemannian manifold $M^n$ when the kinetic energy
is defined by the metric tensor and the form $F$ defines magnetic field.
Also this functional corresponds to the periodic orbits for other
problems of classical mechanics and mathematical physics as it was shown in
\cite{N2,N3,NS}.

When the Lagrangian function
$$
\sqrt{g_{ij} \dot x^i \dot x^j} + A_i \dot x^i
$$
is everywhere positive we obtain a Finsler metric and the periodic
problem can be studied by the methods of the classical Morse
theory. This case is very well understood.
But otherwise we cannot use these methods
because this functional is not bounded from below and also does not
satisfy the Palais-Smale conditions or its analogues.

Its investigation was began by Novikov (\cite{N1,N2,N3}) who,
in particular, introduced a method of throwing
out cycles for overcoming topological difficulties.
He also considered the case when the functional $S$ is multi-valued.
This happens if a magnetic field
is a closed but not exact 2-form.

Some
results about the existence of non-selfintersecting closed extremals on
two-dimensional manifolds were obtained in
\cite{T2,T3} for strong magnetic fields.
A proof of the first result was sketched in \cite{NT}. Later its
generalizations were proved by using of other methods in
\cite{T2,T3}. The methods of these papers strongly use
special properties of the two-dimensional case.

Arnold (\cite{A}) and Kozlov {\cite{Ko}) suggested to use for
this problem methods of symplectic topology (the generalization of
the Poincare theorem which was
proposed by Arnold and was
proved by Conley and Zehnder (\cite{CZ})) . They completed their program  for
positive magnetic field on
the flat two-dimensional torus. Notice that Grinevich and Novikov proved that
the principle of throwing cycles holds in this situation (\cite{GN}).

Symplectic methods  were developed by Ginzburg who obtained
analogous estimates for very strong positive magnetic fields on arbitrary
orientable two-dimen\-sio\-nal manifolds (\cite{G1,G2}).

Notice that these results relate to the case then
the functional $S$ is multi-valued and their proofs also used  specific
properties of the two-dimensional case.

As it was also noticed by Ginzburg there are no periodic trajectories
on two-dimensional Riemann surfaces with genus $\geq 2$ and constant
negative curvature, in the magnetic field which is given by an area
form on this manifold (\cite{G2}). M. Ratner later pointed out that
the paper of Hedlund (\cite{H}) contains a complete description
of such flows on the Lobachevski upperhalf plane for which
the magnetic field is proportional to an area form.


We also notice that in the case of exact magnetic fields negative results
are not known.

In this paper we suggest another approach to this problem and prove the
main theorem which will be formulated below. We are giving here a method of
proving existence of periodic orbits on manifolds with arbitrary
dimension. This method is based on approximation of the functional $S$ by
auxiliary functionals which satisfy the Palais-Smale conditions and getting
closed extremals as limits of sequences of closed extremals of auxiliary
functionals. This works under an additional condition, i.e., in the case when
the ``Ricci curvature'' of this Lagrangian system is positive. One can see
that this approach uses ideas of the seminal paper of Sacks and Uhlenbeck
(\cite{SU}).

Similar phenomena have been studied in the framework of contact forms,
when one has to extremize on Legendre curves for a form $\beta$, the
$L^2$-norm of $\alpha(x)$, where $\alpha$ is another contact form
(see \cite{B}).

But first let us introduce
some definitions.

Let $R^i_{jkl}$ be a Riemann curvature tensor and $\{e_1,\dots,e_n\}$ be
an orthonormal basis in a tangent space $T_xM^n$. Then a Ricci curvature
tensor is defined by the formula
$$
K(u,v) = - \sum^{n}_{\alpha=1} (R(u,e_{\alpha})v,e_{\alpha}).
$$

Now we introduce in the same manner the following tensor
$$
H(w) = \sum^{n}_{\alpha=1} G(e_{\alpha},e_{\alpha},w)
$$
where
$$
G(u,v,w) = (\nabla_uF)_{ij}v^jw^k.
$$

Those definitions, as it can be proved, does not depend on the choice of the
orthonormal bases $\{e_j\}$ in  tangent spaces.

\proclaim{ Main Theorem}
Let $M^n$ be a closed Riemannian manifold endowed with
a Riemannian metric $g_{ij}$ and an exact 2-form $F$ which satisfy the
following condition
$$
\min_{|v|=1}\{f(v) = \frac{K(v,v)}{|v|} - H(v)\}>0.
$$
Then there exists a closed non-one-point extremal of the functional
$$
S(\gamma)=\int_{\gamma}(\sqrt{g_{ij}\dot x^i\dot x^j} + A_i\dot x^i)dt
$$
where $dA=F$.
\endproclaim

This Theorem will be proved in \S 3.

Indeed, the method introduced in the present paper gives more : one almost
can derive the Morse inequalities for everywhere positive functionals and
get the Morse--Novikov inequalities for a  functional which is somewhere
negative
(see Theorem 1).
This needs some additional work.

We would like to notice that the function $f(v)$ for
the case of such Lagrangian systems is the analogue of the Ricci curvature
tensor. So if one can treat a Ricci curvature as an indefinite Riemannian
metric
tensor a {\it Ricci curvature of general homogeneous Lagrangian system}
ought to be treated
as an {\it indefinite Finsler metric}.

\head
\S 2. Functional spaces, auxiliary functionals, and their
analytic properties
\endhead

Let $M^n$ be a differentiable $n$-dimensional closed Riemannian
manifold. We denote by $g_{ij}$ a metric tensor on $M^n$, and denote by
$F$ a closed 2-form on $M^n$.

Let  assume that this form is an exact form, i.e. there exists 1-form
$A_idx^i$ such that $dA=F$:
$$
F_{ij}=\frac{\partial A_j}{\partial x^i}-\frac{\partial A_i}{\partial x^j}.
$$

We denote by $\nabla$ the Levi-Civita derivation, i.e. the covariant
derivation which is compatible with the metric and has vanishing torsion.
The scalar product of tangent vectors $\eta_1,\eta_2 \in T_xM^n$ we denote by
$(\eta_1,\eta_2)$.

Let us denote by $H(S^1,M^n)$ a Hilbert manifold which is formed by maps
$\gamma: [0,1] \rightarrow M^n$ such that
$\gamma(0)=\gamma(1)$. The mapping $\gamma$ is a $H^1$-mapping
(\cite{K1,K2}).
The tangent space $T_{\gamma}H(S^1,M^n)$ of $H(S^1,M^n)$ at
$\gamma$ can be identified with the space of $H^1$-vector fields $\xi (t)$
along the contour $\gamma$.

We denote by $|\dot \gamma|$ a norm of the tangent vector $d\gamma/dt$, and
denote by $\frac{D}{\partial t}$ the covariant derivative along the contour
$\gamma$.
Now we can introduce the Riemann metric $(\xi,\eta)_H$ on
fibers of the tangent bundle $T_{\gamma}H(S^1,M^n)$:
$$
(\xi,\eta)_H = \int^{1}_{0} \{(\xi(t),\eta(t)) + (\frac{D\xi}{\partial t}
(t),\frac{D\eta}{\partial t} (t)) \} dt
$$
where $\xi,\eta \in T_{\gamma}H(S^1,M^n)$.

We see that the group $G$ of orientation-preserved $H^1$-homeomorphisms of
the circle $S^1$ naturally acts on $H(S^1,M^n)$. Let us denote by
$G_0$ its subgroup which is formed by orientation-preserved
$H^1$-homeomorphisms of the interval $[0,1]$. The factor-spaces
$H(S^1,M^n)/G_0$ and $H(S^1,M^n)/G$ we denote by $\Omega(M^n)$ and
$P(M^n)$ respectively. They can be identified with spaces of oriented
non-parameterized $H^1$-contours with based point and unbased point
respectively.

We denote by $H_0(S^1,M^n),\Omega_0(M^n)$, and $P_0(M^n)$ the components, of
these spaces, formed by the contours homotopic to zero.

The scheme of the proof of the main results is the following one. First we
will introduce the auxiliary functionals $S_{\varepsilon,\tau}$ which
becomes the functional
$$
S(\gamma)=\int_{\gamma}(|\dot\gamma| + A_i\dot\gamma^i) dt,
$$
when
$$
\varepsilon,\tau \rightarrow 0.
$$
These functionals will satisfy the Palais-Smale conditions.
After that we will show that for sufficiently small positive values of
$\varepsilon$ and $\tau$ to every nontrivial cycle $H_*(H_0(S^1,M^n),\{S\leq
0\})$ there corresponds a critical point of bounded Morse index of the
functional $S_{\varepsilon,\tau}$.
We will find conditions, on metrics and 2-forms $F$, which will imply
$$
length(\gamma) < const \ index(\gamma)
$$
for critical points of functionals $S_{\varepsilon,\tau}$ for sufficiently
small
values of $\varepsilon$ and $\tau$.
Under these conditions we will prove that there exists a sequence of
critical points $\gamma_{\varepsilon_k,\tau_k}$,
of functionals $S_{\varepsilon_k,\tau_k}$, such that their indices and lengths
of these extremals are uniformly bounded from above,
$\varepsilon_k,\tau_k \rightarrow 0$ as $k \rightarrow\infty$.
In the space
of contours on $M^n$ these extremals tend to  closed non-trivial
extremals of the functional $S$.

Let introduce now the main functionals:
$$
\Phi (\gamma) = \int^{1}_{0} A_j\dot \gamma^j dt,
$$
$$
E_{\theta}:H(S^1,M^n) \rightarrow {\bold R}, 0 \leq \theta \leq 1
$$
where
$$
E_{\theta}(\gamma)=\int^{1}_{0} |\dot \gamma|^{1+\theta} dt.
$$

$E_0$ is the usual functional of length which is $G$-invariant.

The main auxiliary functionals which we will consider are the following
functionals
$$
S_{\varepsilon,\tau} = \varepsilon E_1 + E_{\tau} + \Phi
$$
where $\varepsilon >0, 0 < \tau <1.$

Given a functional $D \in \{E_{\theta}, \Phi\}$ we denote by $\delta
D(\gamma)$ and $\delta^2 D(\gamma)$ its first and second variational
derivatives at $\gamma$ respectively.

The method of computing of these derivatives is very useful
(\cite{K1,K2,M}).
Let us denote by $\exp : TM^n\rightarrow M^n$ the
exponential mapping : $(x,v) \rightarrow exp(x,v)
\in M^n, x \in M^n, v \in T_xM^n$ which corresponds to the Riemannian metric
$g_{ij}$ on $M^n$. By $V_{\gamma,\xi,\eta}(r,s)$ we denote the following
two-dimensional family of contours :

1) $V_{\gamma,\xi,\eta}(r,s):[0,1] \rightarrow M^n,
V_{\gamma,\xi,\eta}(r,s)(0) = V_{\gamma,\xi,\eta}(r,s)(1)$ ;

2) $V_{\gamma,\xi,\eta}(r,s)(t) = \exp(\gamma(t), r\xi(t) + s\eta(t)),$

\noindent
where $\xi, \eta \in T_{\gamma}H(S^1,M^n), r,s \in {\bold R}.$

Since any closed curve which is sufficiently close to
$\gamma$ is represented in the form $V_{\gamma,\xi,0}(1,0)(t)$
and since the exponential mapping determine local charts in the Hilbert
manifold
$H(S^1,M^n)$, we have the following  lemma.

\proclaim{ Lemma 1}
$$
\delta S(\gamma)(\xi) = (\partial S(V_{\gamma,\xi,0}(r,0))/ \partial r)
\vert_{r=0},
$$
$$
\delta^2S(\gamma)(\xi,\eta) = (\partial^2S(V_{\gamma,\xi,\eta}(r,s))/ \partial
r \partial s)\vert_{r=s=0}.
$$
\endproclaim

Now we are able to compute variational derivatives of the functionals
$S_{\varepsilon,\tau}$.

\proclaim{ Lemma 2}
Let $\theta >0$. Then $E_{\theta} \in C^1(H(S^1,M^n),
{\bold R})$ and
$$
\delta E_{\theta}(\gamma)(\xi) = (1+\theta) \int^{1}_{0}
(\frac{D\xi}{\partial t},\dot \gamma) |\dot \gamma|^{\theta-1} dt.
$$
If also $\gamma \in C^2(S^1,M^n)$ and the vector-field $\xi(t)$ is of class
$C^1$ then the formula
$$
 \delta E_{\theta}(\gamma)(\xi) = -(1+\theta) \{\int^{1}_{0}
(\xi,\frac{D\dot \gamma}{\partial t}) |\dot \gamma|^{\theta-1} dt +
 $$
$$
+ (\theta-1)\int^{1}_{0} (\xi,\dot \gamma)(\frac{D\dot \gamma}{\partial
t},\dot \gamma)|\dot \gamma|^{\theta-3}dt \}
$$
holds.
\endproclaim

{\bf Proof of the Lemma 2.}

For the sake of simplicity, let us denote $V_{\gamma,\xi,0}$ by $V_{\gamma}$.

Since the covariant derivation is compatible with the Riemannian metric
\newline
($\nabla g_{ij} \equiv 0$),
$$
\frac{D}{\partial t}(\alpha,\beta)=(\frac{D\alpha}{\partial t},\beta) +
(\alpha,\frac{D\beta}{\partial t})
$$
for every pair of $H^1$-vector-fields $\alpha$ and $\beta$ along $\gamma$.
In the sequel we will use this formula for
vector-fields along curves of the form
$\ r \rightarrow exp(\gamma(t),r\xi(t))$.
Since $V_{\gamma,\xi,0}$ can be considered to be a surface parametrized
by $t$ and $r$, the Lie bracket
$[\frac{\partial}{\partial t},\xi]$ is zero for a suitable parametrization.
Therefore, the covariant derivation has vanishing torsion and we derive
$$
\nabla_{\xi}\frac{\partial V_{\gamma}}{\partial t} = \frac{D}{\partial
r}\frac{\partial V_{\gamma}}{\partial t} = \frac{D\xi}{\partial t}.
$$
Now by using of Lemma 1 we obtain
$$
\delta E_{\theta}(\gamma)(\xi) = \frac{\partial}{\partial r}(\int^{1}_{0}
(\dot V_{\gamma},\dot V_{\gamma})^{(\theta+1)/2} dt) \vert_{r=0} \ =
 $$
$$
(1+\theta) \int^{1}_{0} ((\nabla_{\xi}\dot V_{\gamma}, \dot
\gamma)\vert_{r=0})|\dot \gamma|^{\theta-1} dt =
(1+\theta)\int^{1}_{0}(\frac{D\xi}{\partial t},\dot \gamma)|\dot
\gamma|^{\theta-1} dt,
$$
where we denote by $\dot V_{\gamma}$ the covariant derivative $\frac{\partial
V_{\gamma}}{\partial t}$.

$\gamma \in C^2(S^1,M^n)$ and the vector-field $\xi$ is of class $C^1$.
We integrate by part and derive
$$
(1+\theta)\int^{1}_{0}(\frac{D\xi}{\partial t},\dot \gamma)|\dot
\gamma|^{\theta-1} dt =
$$
$$
(1+\theta)\int^{1}_{0} (\frac{D}{\partial t} (\xi,\dot \gamma) -
(\xi,\frac{D}{\partial t}\dot \gamma))|\dot \gamma|^{\theta-1} dt =
$$
$$
-(1+\theta)\int^{1}_{0}(\xi,\frac{D}{\partial t}\dot \gamma)|\dot
\gamma|^{\theta-1} dt - (1+\theta)\int^{1}_{0}(\xi,\dot
\gamma)(\frac{D}{\partial t}|\dot \gamma|^{\theta-1})dt.
$$
We have
$$
\frac{D}{\partial t}|\dot \gamma|^{\theta-1} =
(\theta-1)(\frac{D}{\partial t}\dot \gamma,\dot \gamma)|\dot
\gamma|^{\theta-3}.
$$

The claim about the smoothness of $E_{\theta}$ follows from
the formula for $\delta E_{\theta}$ immediately.

Lemma 2 is established.

\proclaim{ Lemma 3}
$$
\delta \Phi (\gamma)(\xi) = \int^{1}_{0} F_{ij}\xi^i\dot \gamma^j dt.
$$
\endproclaim

{\bf Proof of Lemma 3.}

Let us divide the contour $\gamma$ into arcs $\gamma_{\alpha}, 1 \leq \alpha
\leq k,$ such that in a neighborhood of every arc $\gamma_{\alpha}$ it is
possible to introduce local coordinates and make sense of expressions of the
form $\partial \xi^j/\partial t$ which are not covariant.

$$
\delta \Phi (\gamma)(\xi) = \frac{d}{dr}(\int^{1}_{0} A_i\dot
V^i_{\gamma}(r,0))\vert_{r=0} =$$
$$
\int^{1}_{0} (\nabla_{\xi}A_i)\dot \gamma^i dt + \int^{1}_{0}
A_i(\nabla_{\xi} \dot V^i_{\gamma}(r,0))\vert_{r=0} dt =
$$
$$
\sum^k_{\alpha=1}\{\int_{\gamma_{\alpha}}(\partial_jA_i -
\Gamma^k_{ji}A_k)\xi^j\dot\gamma^i dt + \int_{\gamma_{\alpha}}
A_i\frac{D\xi}{\partial t}^i dt\} =
$$
$$
\sum^k_{\alpha=1}\{\int_{\gamma_{\alpha}}(\partial_jA_i -
\Gamma^k_{ji}A_k)\xi^j\dot\gamma^i dt + \int_{\gamma_{\alpha}}
A_i(\frac{\partial \xi^i}{\partial t}+\Gamma^I_{jk}\xi^j\dot\gamma^k)dt \}=
$$
$$
\sum^k_{\alpha=1}\int_{\gamma_{\alpha}} (\partial_jA_i\xi^j\dot\gamma^i +
A_i\frac{\partial \xi^i}{\partial t}) dt
$$
where $\Gamma^i_{jk}$ are the Christoffel symbols which correspond to the
connectedness $\nabla$.

Using integration by parts we obtain
$$
\delta \Phi(\gamma)(\xi) = \sum^k_{\gamma{\alpha}} \int_{\gamma_{\alpha}}
(\partial_iA_j - \partial_jA_i)\xi^i\dot\gamma^j dt = \int^{1}_{0}
F_{ij}\xi^i\dot\gamma^j dt.
$$

Lemma 3 established.

To compute the Morse indices of the extremals of functionals
$S_{\varepsilon,\tau}$ it is sufficient to find out the
second derivatives of these
functionals at the critical points.

We start with the following fact.
If $\gamma$ is a critical point of the functional
$S_{\varepsilon,\tau}$ then the motion of a point along this contour is
determined by the Euler-Lagrange equations
with the Lagrange function
$$
L_{\varepsilon,\tau}(x,\dot x) = \varepsilon |\dot x|^2 + |\dot x|^{1+\tau} +
A_i(x)\dot x^i.
$$

The following lemma holds.

\proclaim{Lemma 4}
 If $x(t)$ is a trajectory of a system with the
Lagrange function
$L_{\varepsilon,\tau}$, where $\varepsilon,\tau >0$, then $|\dot x| =
const.$
\endproclaim

{\bf Proof of Lemma 4.}

By the Euler-Lagrange equations,
the function
$I=\dot x^i \frac{\partial L}{\partial \dot x^i} -L$
is constant along the trajectories. In our case
$I=\varepsilon|\dot x|^2+\tau|\dot x|^{1+\tau},$ and thus we obtain
that $|\dot x|$ is constant along the trajectories.

Lemma 4 follows.

The following Lemma now immediately follows from Lemmas 2, 3, and 4
that give this form for the Euler-Lagrange equation.

\proclaim{Lemma 5}
The extremals of the functional $S_{\varepsilon,\tau}$ satisfy the following
equation
$$
\big(\frac{D \dot \gamma}{\partial t}\big)^j =
\frac{\partial}{\partial t} \dot \gamma + \Gamma^j_{ik}
\dot\gamma^i \dot \gamma^k =
\frac{g^{jk}F_{ki}\dot \gamma^i}{2\varepsilon +
(1+\tau)|\dot\gamma|^{\tau-1}}.
\eqno{(2.1)}
$$
\endproclaim

Let us now recall the definition of the Riemann curvature tensor:
$$
R^i_{jkl}\xi^j\eta^k\zeta^l = [R(\xi,\eta)\zeta]^i =
(\nabla_{\xi}\nabla_{\eta}\zeta - \nabla_{\eta}\nabla_{\xi}\zeta -
\nabla_{[\xi,\eta]}\zeta)^i,
$$
where $\xi,\eta$ and $\zeta$ are smooth vector-fields
on $M^n$.

\proclaim{Lemma 6}
Let
$\gamma$ be a nontrivial critical point of the functional $S_{\epsilon,\tau}$.
Then at the point $\gamma \in H(S^1,M^n)$ the second variational derivative
of the functional $S_{\varepsilon,\tau}$ exists and
$$
\delta^2S_{\varepsilon,\tau}(\gamma)(\xi,\eta) =
$$
$$
2\varepsilon\int_0^1 \big(\big(\frac{D\xi}{\partial t},\frac{D\eta}{\partial t}
\big) + (R(\eta,\dot\gamma)\xi,\dot\gamma)\big) dt +
$$
$$
(1+\tau)\int^{1}_{0}\{(\tau-1)|\dot \gamma|^{\tau-3}(\frac{D\xi}
{\partial t},\dot \gamma)(\frac{D\eta}{\partial t},\dot \gamma) +
\eqno{(2.2)}
$$
$$
|\dot \gamma|^{\tau-1}((\frac{D\xi}{\partial t},\frac{D\eta}{\partial
t})+(R(\eta,\dot \gamma)\xi,\dot \gamma))\}dt +
$$
$$
\int_0^1 \{(\nabla_{\eta}F)_{jk}\xi^j \dot\gamma^k + F_{jk}\xi^j
\big(\frac{D \eta}{\partial t}\big)^k\}dt.
$$
\endproclaim

{\bf Proof of Lemma 6.} By Lemmas 1 and 2, we have
$$
\delta^2E_{\theta}(\gamma)(\xi,\eta) = (1+\theta)\int^{1}_{0}
\frac{D}{\partial s}((\frac{D\xi}{\partial t},\dot V_{\gamma})|\dot
V_{\gamma}|^{\theta-1})\vert_{s=r=0} dt =
$$
$$
(1+\theta)\int^{1}_{0} \{(\theta-1)|\dot
\gamma|^{\theta-3}(\frac{D}{\partial s}\dot V_{\gamma},\dot
\gamma)(\frac{D\xi}{\partial t},\dot \gamma) + $$
$$|\dot \gamma|^{\theta-1}((\frac{D}{\partial s}\frac{D}{\partial t}\xi,\dot
\gamma) +(\frac{D\xi}{\partial t},\frac{D}{\partial s}\dot V_{\gamma}))\} dt
\vert_{s=r=0}=
$$
$$
(1+\theta)\int^{1}_{0} \{(\theta-1)|\dot
\gamma|^{\theta-3}(\frac{D\xi}{\partial t},\dot
\gamma)(\frac{D\eta}{\partial t},\dot\gamma) +
\eqno{(2.3)}
$$
$$
|\dot\gamma|^{\theta-1}(\frac{D}{\partial t}\frac{D}{\partial s}\xi +
R(\eta,\dot\gamma)\xi,\dot\gamma)+|\dot\gamma|^{\theta-1}
(\frac{D\xi}{\partial t},\frac{D\eta}{\partial t}) \}dt.
$$

Here we use the following fact : if we introduce the coordinates
$s$ and $t$ on a surface immersed into $M^n$ (in our case $(s,t)
\rightarrow \exp(\gamma(t),s\eta(t))$) then, by definition of
the Riemannian curvature, we have
$$
\frac{D}{\partial s}\frac{D}{\partial t}\xi - \frac{D}{\partial
t}\frac{D}{\partial s}\xi = R(\eta,\dot \gamma)\xi.
$$

By Lemmas 1 and 3, we also have
$$
\delta^2\Phi(\gamma)(\xi,\eta) = \int^{1}_{0}\frac{D}{\partial
s}(F_{jk}\xi^j\dot (V_{\gamma})^k)\vert_{s=r=0} dt=
$$
$$
\int^{1}_{0}\{ (\nabla_{\eta}F)_{jk}\xi^j\dot \gamma^k +
F_{jk}(\nabla_{\eta}\xi)^j\dot \gamma^k + F_{jk}\xi^j(\nabla_{\eta}\dot
\gamma)^k\}dt.
\eqno{(2.4)}
$$

But
$$
\nabla_{\eta}\dot \gamma = \frac{D}{\partial t}\eta.
$$

We  obtain the formula for the second derivative of $S_{\varepsilon,\tau}$
by adding the  formulas for second derivatives of
$\varepsilon E_1$, $E_{\tau}$ and $\Phi$. We do not give here this formula  but
show how it reduces to (2.2).

Let us consider terms in this formula that contain
$\nabla_{\eta}\xi = \frac{D}{\partial s}\xi$:
$$
J = \int_0^1 (2\varepsilon + (1+\tau)|\dot\gamma|^{\tau-1})
\big(\frac{D}{\partial t}\nabla_{\eta}\xi,\dot\gamma\big) dt +
\int_0^1 F_{jk} (\nabla_{\eta}\xi)^j \dot\gamma^k dt.
\eqno{(2.5)}
$$
Integrating the first summand by parts and by using (2.1) we
derive that $J=0$.

The sum of the other terms gives (2.2).

Lemma 6 follows.

The following lemma is a consequence of Lemmas 4 and 6.

\proclaim{Lemma 7}
Let $\gamma$ be a nontrivial critical point of the functional
$S_{\varepsilon,\tau}$.  Then
there exists a neighborhood of $\gamma$ such that
in this neighborhood the second derivative of the functional
$S_{\varepsilon,\tau}$ exists and is defined by formula
(2.2) on this whole neighborhood.
\endproclaim

It is easy to see that the neighborhood formed by contours which are close
to $\gamma$ and which have a never vanishing tangent vector
satisfy to requirements of this lemma.

\head
\S 3. Proof of the main theorems
\endhead

{\bf A)} First let us prove the following fact.

\proclaim{Lemma 8}
The subset of $H(S^1,M^n)$ formed by one-point-contours
is a manifold of local minima of the functional
$S_{\varepsilon,\tau}$, $0 \leq \tau <1$.
\endproclaim

{\bf Proof of Lemma 8.} Using the Holder inequality,
$$
(\int^{1}_{0}f(t)dt)^{1+\tau} \leq \int^{1}_{0}f(t)^{1+\tau}dt
$$
for  $f \in L_2[0,1], f \geq 0$.
Substituting  $|\dot \gamma|$ for $f(t)$ and taking into account that by
well-known results from the theory of isoperimetric inequalities
the following inequality
$$
\int_{\Pi}F \leq c \ length(\gamma)^2,\eqno{(3.1)}
$$
holds for sufficiently small contours $\gamma$, where $\Pi$ is a
sufficiently small surface having $\gamma$  as boundary,
$c$ is a positive constant which depends on the metric
$g_{ij}$ and the 2-form $F$
we obtain the inequality
$$
S_{\varepsilon,\tau} \geq length(\gamma)^{1+\tau} - c \
length(\gamma)^2,\eqno{(3.2)}
$$
which holds for all sufficiently small contours from
$H_0(S^1,M^n)$ since every such contour can be transformed into a smooth
contour by a perturbation as small as we may wish.
The claim of Lemma 8 follows from the equality (3.2).

Lemma 8 is proved.

\proclaim{Definition}
 A nontrivial cycle $\alpha \in H_{*}(M^n)$ can be
thrown out into the domain $\{S_{\varepsilon,\tau}<0\}$
if and only if there exists a
continuous mapping
$$
\phi: K \times [0,1] \rightarrow H_0(S^1,M^n)
$$
such that

1) $K$ is a polyhedron and $[\phi(K \times 1)] = \alpha$, i.e. an image of
$K \times 0$ realizes a cycle $\alpha$ ;

2) $\phi(K \times 0) \subset \{S_{\varepsilon,\tau}<0\}$.
\endproclaim

The following lemma holds.

\proclaim{Lemma 9}
 If $\tau<1$ and a cycle $\alpha \in H_k(M^n) \setminus
\{0\}$ can be thrown out into $\{S_{\varepsilon,\tau}<0\}$ then

1) there exists a positive constant $c_{\varepsilon,\tau}(\alpha)$
(the ``critical value'' of the functional $S_{\varepsilon,\tau}$)
defined by the formula
$$
c_{\varepsilon,\tau}(\alpha) = \inf_{\phi \in A(\alpha)}
\max \{S_{\varepsilon,\tau}(\gamma) :
\gamma \in Im \phi \} >0
$$
where $A(\alpha)$ is the set of all possible throwings out of a cycle
$\alpha$;

2) $H_{k+1}(H_0(S^1,M^n),\{S_{\varepsilon,\tau} \leq 0\}) \neq 0$;

3) if $D_k$ is the linear span of $k$-cycles (with coefficients
from the field $E$) which can be thrown out
into $\{S_{\varepsilon,\tau} <0\}$, then
$$
\dim_E H_{k+1}(H_0(S^1,M^n),\{S_{\varepsilon,\tau} \leq 0\};E) \geq
\dim_E D_k.
$$
\endproclaim

For proof of this lemma see \cite{T1,T2}.
We also can use the results of these papers in the sequel.

\proclaim{Lemma 10}
 Let there exist a contour $\gamma \in \Omega_0(M^n)$
such that $(E_0+\Phi)(\gamma)<0$. Then there exist positive constants
$\varepsilon_0$ and $\tau_0$ such that all the manifold $M^n$
is thrown out into $S_{\varepsilon,\tau}<0$ for $\varepsilon < \varepsilon_0,
\tau<\tau_0$.
\endproclaim

{\bf Proof of Lemma 10.}

By the construction of the throwing out given in \cite{T1},
the contours from
$\phi(M^n \times [0,1])$ are piecewise-smooth and
we can assume that they lie in $H_0(S^1,M^n)$.

Since the subset $\phi (M^n \times 0)$ is the continuous image of the compact
space, there exist sufficiently small positive constants
$\varepsilon_0$ and $\tau_0$ such that
for $\varepsilon<\varepsilon_0, \tau<\tau_0$ and for any contour $\gamma \in
\phi(M^n \times 0)$ the following inequality
$$
(\varepsilon E_1 + E_{\tau} + \Phi)(\gamma) <0
$$
holds.

Lemma 10 is proved.

\proclaim{ Lemma 11}
Let the conditions of Lemma 10 hold. Assume  that
the manifold $M^n$ can be thrown out into
$\{S_{\varepsilon,\tau}<0\}$ for $0<\varepsilon<\varepsilon_0, 0<\tau,\tau_0$.
Then we can define map which corresponds to non-zero cycles of $H_*(M^n)$
critical values of the functional $S_{\epsilon,\tau}$:
$$
\alpha \in H_{*}(M^n) \setminus \{0\} \rightarrow
c_{\varepsilon,\tau}(\alpha).
$$
This mapping is semicontinuous above with respect to
variables $\varepsilon$ and $\tau$.
\endproclaim

{\bf Proof of Lemma 11.} It is sufficient to prove
that for any constants
$\varepsilon,\tau$ such that $\varepsilon \in (0,\varepsilon_0), \tau \in
(0,\tau_0)$ and for any $\omega >0$ there exist  positive constants
$\lambda$ and  $\mu$ such that
$$
c_{\varepsilon',\tau'}(\alpha) - c_{\varepsilon,\tau}(\alpha) < \omega
$$
for $|\varepsilon-\varepsilon'| < \lambda$ and $|\tau-\tau'| < \mu$.

By the definitions of the throwing out of cycles and critical values
$c_{\varepsilon,\tau}(\alpha)$ corresponding to them,
there exists the continuous mapping
$$
\phi : K \times [0,1] \rightarrow H_0(S^1,M^n),
$$
such that

1) $[\phi(K \times 1] = \alpha$ ;

2) $\phi(K \times 1) \subset \{S_{\varepsilon,\tau}<0\}$;

3) $c_{\varepsilon,\tau} \leq \max \{S_{\varepsilon,\tau}(\gamma) :
 \gamma \in Im
\phi\} \leq c_{\varepsilon,\tau} + \omega/2$.

Since the subspaces $Im \phi$ and $\phi(K\times 0)$ are  continuous
images of compact spaces, it is evident that such constants
$\lambda$ and $\mu$ exist.

Lemma 11 follows.

\proclaim{Lemma 12}
Let $\varepsilon$ and $\tau$ be sufficiently small.
Then there exist a positive constant $J_0$ such that the
lengths of
extremals of the functionals $S_{\epsilon,\tau}$ which are
not points
are not less than $J_0$.
\endproclaim

{\bf Proof of Lemma 12.}

By Lemma 5, extremals of the functional $S_{\varepsilon,\tau}$
satisfy the Euler-Lagrange equation
$$
\frac{\partial}{\partial t} \dot \gamma^j +
\Gamma^j_{ik}\dot\gamma^i\dot\gamma^k =
g^{jk}F_{ki}\dot\gamma^i/
(2\varepsilon + (1+\tau)|\dot\gamma|^{\tau-1}).\eqno{(3.3)}
$$

By Lemma 4, these extremals are naturally parameterized and
$$
length(\gamma) = |\dot\gamma| = const.
$$

Let $\lambda$ be sufficiently small and let an extremal $\gamma$ lie in a
small neighborhood of the contour $\gamma(0)$ endowed with local coordinates
$(x^1,\dots,x^n)$ with
$\dot\gamma = (\lambda,0,\dots,0), \lambda = |\dot\gamma|$.
By (3.3),
$$
\dot\gamma(t) = \dot\gamma(0) + o(\lambda^{3/2}),
$$
for $\tau<1/2$, and, thus, if $\lambda$ is sufficiently small then
$\dot\gamma^1$ is nowhere equal to zero and an extremal $\gamma$ can not be
closed.

Lemma 12 follows.

{\bf B) Validity of the Palais-Smale conditions for the functionals
$S_{\varepsilon,\tau}$.}

\proclaim{(PS) Condition (\cite{PS})}
 A $C^1$-functional $S$ on a
Hilbert manifold $X$ satisfies the (PS) condition if
and only if for every subset $W
\subset X$ on which the functional is bounded from above and the norms of
its first derivatives $\|S\|$ are not bounded from below there exist
a point $x$ which lies in the closure of $W$
satisfying $\delta S(x)=0$.
\endproclaim

First let us mention very simple fact.

\proclaim{Lemma 13}
 $E^2_0(\gamma) \leq E_1(\gamma)$
for $\gamma \in H_0((S^1,M^n)$.
\endproclaim

This lemma is an immediate consequence
of the well-known inequality
$$
\int^{1}_{0} fg dx \leq \sqrt{\int^{1}_{0} f^2 dx}
\sqrt{\int^{1}_{0} g^2 dx}
$$
where $f,g \in L_2[0,1]$.

\proclaim{Lemma 14}
There exists a positive constant $T$ such that
$$
E_1(\gamma) \leq (T+ \sqrt{T^2 + 4\varepsilon S_{\varepsilon,\tau}})^2/
4\varepsilon^2.
$$
\endproclaim

{\bf Proof of Lemma 14.}
Since the manifold $M^n$ is compact, there exists a constant
$T$ such that
$$
|A_jv^j| \leq T|v|
$$
for any $x \in M^n, v \in T_xM^n$.

Since the functional $E_{\tau}$ is non-negative,
$$
\varepsilon E_1 + \Phi \leq \varepsilon E_1 + E_{\tau} + \Phi =
S_{\varepsilon,\tau}.
$$

Now it is easy to see that
$$
\varepsilon E_1 - TE_0 \leq \varepsilon E_1 - |\Phi| \leq S_{\varepsilon,\tau},
$$
and , by Lemma 13, we obtain
$$
\varepsilon E_1 - TE_0 \leq \varepsilon E_1 - T\sqrt{E_1} \leq S_{\varepsilon,
\tau}.
$$

The claim of Lemma 14 follows from the last inequality.

Lemma 14 follows.

\proclaim{Lemma 15}
Let $\varepsilon>0, 0<\tau<1$. Then the functional
$$
S_{\varepsilon,\tau} : H_0(S^1,M^n) \rightarrow M^n
$$
satisfy (PS) conditions.
\endproclaim

{\bf Proof of Lemma 15.}
It is sufficient to show that for every sequence
$\{\overline{\gamma_j}\}\subset H_0(S^1,M^n)$ such that

1) there exists constant $K_0$ such that $S_{\varepsilon,\tau}(\gamma_j) < K_0
< +\infty$ for $j=1,2,\dots;$

2) $\|\delta S_{\varepsilon,\tau}(\overline{\gamma_j})\| \rightarrow 0$ as
$j \rightarrow +\infty,$
\newline
there exists a subsequence $\{\gamma_j\}$ which converges to a critical point
of the functional $S_{\varepsilon,\tau}$.

Here we denote by $\|\delta S_{\varepsilon,\tau}\|$ a norm of the variational
differential:
$$
\|\delta S_{\varepsilon,\tau}(\gamma)\| = \sup \{|\delta S_{\varepsilon,\tau}
(\nu)|: \nu\in T_{\gamma}H_0(S^1,M^n), \|\nu\| =1\}.
$$

Let us consider such a sequence $\{\overline{\gamma_j}\}$.

By Lemma 14, there exists a constant $K_1$ such that
$$
E_1(\overline{\gamma_j}) \leq K_1
$$
for $j=1,2,\dots.$

Let us denote by $d=d_M$ the distance function, on $M^n$, generated by the
Riemannian metric
$g_{ij}$. By the Cauchy-Schwarz inequality,
$$
d(\overline{\gamma}_i(t_0),\overline{\gamma}_i(t_1)) \leq
\sqrt{2|t_0-t_1|E_1(\overline{\gamma}_i)} \leq \sqrt{2K_1|t_0-t_1|}.
$$

The last inequality implies the sequence
$\{\overline{\gamma}_j\}$ is equicontinuous and, since the sequence
$\{\overline{\gamma}_j(t_0)\}$ is relatively compact for any $t_0$, then, by
the Arzela-Ascoli theorem, there exists a subsequence
$\{\gamma_k\} \subset \{\overline{\gamma}_j\}$
which converges to some contour $\gamma^{*} \in C^0(S^1,M^n)$ in the
space $C^0(S^1,M^n).$

It remains to prove that this $\gamma^{*} \in H_(S^1,M^n)$
satisfies
$$
\|\delta S_{\varepsilon,\tau}(\gamma^{*})\|=0.
$$

Let us consider small perturbation $\omega$ of the contour $\gamma^{*}$
such that

1) $\omega \in C^{\infty}(S^1,M^n)$;

2) $\gamma_j(t) = \exp(\omega(t),\xi_j(t))$ where
$\xi_j \in T_{\omega}H_0(S^1,M^n)$, $|\xi_j(t)|_M \leq K_2$ for any $t \in
[0,1]$, $j \geq j_0;$

3) the exponential mapping $exp:
T_{\omega(t)}M^n \rightarrow M^n$ is a diffeomorphism onto its image
on balls $|\eta|_M \leq 10K_2 < inj(M^n)$ for any $t \in [0,1]$
and on every such ball the mappings
$$
\frac{\partial\ \exp(\omega(t),\eta)}{\partial x}: T_{\omega(t)}M^n
\rightarrow T_{\exp(\omega(t),\eta)}M^n,
$$
and
$$
\frac{\partial \ \exp(\omega(t),\eta)}{\partial v}: T_{\eta}T_{\omega(t)}M^n
\rightarrow T_{\exp(\omega(t),\eta)}M^n
$$
are invertible.

Without loss of generality, we assume that $j_0=1$.

Since $\gamma_j \rightarrow \gamma^{*}$ in $C^0(S^1,M^n)$,  $\xi_j
\rightarrow \xi^{*}$  in the space of $C^0$-vector-fields:
$$
\rho_j = \max |\xi_j(t)-\xi^{*}(t)|_M \rightarrow 0
$$
as $j \rightarrow + \infty,$ and we derive now that
$$
\int^{1}_{0}|\xi_j(t) - \xi^{*}(t)|^2 dt \rightarrow 0
$$
as $j \rightarrow +\infty.$

It remains to prove that
$$
\int^{1}_{0} |\nabla_{\omega}\xi_i(t) - \nabla_{\omega}\xi_j(t)|^2 dt
\rightarrow 0
$$
as $i,j \rightarrow +\infty$ where $\nabla_{\omega}$ is a covariant
derivation along the contour $\omega.$

For the sake of simplicity, let us put
$$
A_{1,j}(t) = \frac{\partial\ \exp(\omega(t),\xi_j(t))}{\partial x},
\ \ \ A_{2,j}(t) = \frac{\partial\ \exp(\omega(t),\xi_j(t))}{\partial v},
$$
$$
A_1(t) = \frac{\partial\ \exp(\omega(t),\xi^{*}(t))}{\partial x},
\ \ \ A_2(t) = \frac{\partial\ \exp(\omega(t),\xi^{*}(t))}{\partial v}.
$$

Let us represent $\dot \gamma_j(t)$ in the following form
$$
\dot \gamma_j(t) = A_{1,j}(t)\dot \omega(t) +
A_{2,j}(t)\nabla_{\omega}\xi_j(t).
$$
It follows from this representation that
$$
\nabla_{\omega}\xi_j = A^{-1}_{2,j}(\dot \gamma_j - A_{1,j}\dot\omega).
$$
Since the right-hand side of this equality is square integrable
(i.e. $\in L_2$), since $S_{\varepsilon,\tau}(\gamma_j)$ are
uniformly bounded ($<K_0$), using Lemma 14 and the representation
of $\gamma_j$ in the exponential form ($\gamma_j = \exp(\omega,\xi_j)$), we
conclude that
$$
\|\xi_j\|_H=\int^{1}_{0}(|\xi_j(t)|^2+|\nabla_{\omega}\xi_j(t)|^2)dt \leq
K_3 < +\infty
$$
where $K_3$ is some constant.

Let us introduce the vector-fields $u_{ij}$ along $\gamma_i$ defining them
by the following formula
$$
\gamma_j(t) = \exp(\gamma_i(t),u_{ij}(t)).
$$

We obtain that
$$
\nabla_iu_{ij}=\left(\frac{\partial\ \exp(\gamma_i,u_{ij})}{\partial v}
\right)^{-1}(A_{2,j}\nabla_{\omega}\xi_j -
\frac{\partial \ \exp(\gamma_i,u_{ij})}{\partial x}A_{2,i}
\nabla_{\omega}\xi_i) + v_{ij}
$$
where
$$
v_{ij}=\left(\frac{\partial \ \exp(\gamma_i,u_{ij})}{\partial v}
\right)^{-1}(A_{1,j} - \frac{\partial \ \exp(\gamma_i,u_{ij})}{\partial x}
A_{1,i})\dot \omega
$$
and $\nabla_i$ is the covariant derivation along $\gamma_i$.

Since by the construction the contours $\gamma_i$ lie in the sufficiently small
neighborhood of $\omega$, we assume that the norms of mappings
$$
\frac{\partial \ \exp(\gamma_i,u_{ij})}{\partial x},\ \ \
\frac{\partial \ \exp(\gamma_i,u_{ij})}{\partial v}
$$
and others which are inverse to them are bounded by some constant.
Let us notice also that
$$
v_{ij}=O(\rho_{ij})=\max |\xi_i(t)-\xi_j(t)|,
$$
and, thus, there exists constant $K_4$ such that
$$
\int^{1}_{0}(|u_{ij}(t)|^2+|\nabla_iu_{ij}(t)|^2)dt \leq K_4 < +\infty.
$$

Let us find out an explicit formula for
$$
w_{ij} =\delta S_{\varepsilon,\tau}(\gamma_i)(u_{ij}) + \delta
S_{\varepsilon,\tau}(\gamma_j)(u_{ji}).
$$

By Lemmas 2 and 3,
$$
\delta S_{\varepsilon,\tau}(\gamma_i)(u_{ij}) = \int^{1}_{0}
\{(2\varepsilon+(1+\tau)|A_{1,i}\dot \omega +
A_{2,i}\nabla_{\omega}\xi_i|^{\tau-1}) \times
$$
$$
(A_{1,i}\dot \omega + A_{2,i}\nabla_{\omega}\xi_i, \nabla_iu_{ij}) +
F_{kl}u^k_{ij}\dot\gamma^l_i\}dt.
$$

Since $|u_{ij}|=O(\rho_{ij})$, we assume that
$$
\delta S_{\varepsilon,\tau}(\gamma_i)(u_{ij} = \int^{1}_{0}
(2\varepsilon+(1+\tau)|A_1\dot\omega +A_2\nabla_{\omega}\xi_i|^{\tau-1})
\times
$$
$$
(A_1\dot\omega + A_2\nabla_{\omega}\xi_i,
A_2(\nabla_{\omega}\xi_j-\nabla_{\omega}\xi_i)) dt +
O(\rho_{ij}).
\eqno{(3.4)}
$$

Let us prove now the following auxiliary proposition.

\proclaim{ Proposition 1}
 If $u,v \in R^n \setminus \{0\}$ and $0<\tau<1$
then
$$
(|u|^2 -(u,v))|u|^{\tau-1} + (|v|^2-(u,v))|v|^{\tau-1} \geq 0.
\eqno{(3.5)}
$$
\endproclaim

{\bf Proof of Proposition 1.} If $(u,v)\leq 0$ then the claim of proposition
is evident. We are left with the case when
$(u,v)=\beta |u||v|, 0\leq\beta \leq 1.$ Without loss of generality, we assume
that $|u|/|v|=k \geq 1.$ Let us divide the left-hand side of equality (3.5)
by $|v|^{\tau+1}$ and reduce our problem to proving  the inequality
$$
k^{1+\tau} - k^{\tau}\beta \geq k\beta -1.
$$
But $k^{\tau}(k-\beta) \geq (k-\beta)$, and, thus, it is sufficiently to
prove that $(k-\beta)\geq (k\beta-1).$ The last inequality is equivalent
to the following
$$
k(1-\beta)\geq (\beta-1)
$$
which holds since $k\geq 1$  $0\leq\beta\leq 1.$

Proposition 1 follows.

Let us apply inequality (3.5) to vectors $u(t)=A_1(t)\dot\omega(t) +
A_2(t)\nabla_{\omega}\xi_i(t)$ and $v(t)=A_1(t)\dot\omega(t) +
A_2(t)\nabla_{\omega}\xi_j(t)$ from $T_{\gamma^{*}(t)}M^n$. We
obtain that
$$
D_{1,ij}=\int^{1}_{0}(|u|^{\tau+1}(u,v-u)+|v|^{\tau-1}(v,u-v))dt \leq 0.
$$

By (3.4),
$$
w_{ij}= D_{1,ij}+D_{2,ij}+O(\rho_{ij})
$$
where
$$
D_{2,ij}=\int^{1}_{0} 2\varepsilon\{(A_1\dot\omega +
A_2\nabla{\omega}\xi_i,A_2\nabla_{\omega}(\xi_j-\xi_i) +
$$
$$
(A_1\dot\omega + A_2\nabla_{\omega}\xi_j,
A_2\nabla_{\omega}(\xi_i-\xi_j))\} dt.
$$

{}From the expression for $D_{2,ij}$ we conclude that
$$
-D_{2,ij} = 2\varepsilon\int^{1}_{0}|A_2\nabla_{\omega}(\xi_i-\xi_j)|^2 dt.
$$

Now it follows from (3.4) that
$$
2\varepsilon\int^{1}_{0}|A_2\nabla_{\omega}(\xi_i-\xi_j)|^2 dt \leq -w_{ij} +
O(\rho_{ij}) \leq
$$
$$
R_{ij} = K_4(\|\delta S_{\varepsilon,\tau}(\gamma_i)\| + \|\delta
S_{\varepsilon,\tau}(\gamma_j)\|) + O(\rho_{ij}).
$$

But $R_{ij} \rightarrow 0$ as $i,j \rightarrow +\infty$, and, thus,
$$
\int^{1}_{0}|A_2\nabla_{\omega}(\xi_i-\xi_j)|^2 dt \rightarrow 0
$$
as $i,j \rightarrow +\infty$,
which in turn leads to
$$
\int^{1}_{0}|\nabla_{\omega}\xi_i - \nabla_{\omega}\xi_j|^2 dt
 \rightarrow 0
$$
as $i,j \rightarrow +\infty.$

Using that  and the fact that
$$
\int^{1}_{0}|\xi_i-\xi_j|^2 dt \rightarrow 0, \ \ i,j \rightarrow +\infty
$$
we derive that the sequence $\{\xi_i\}$ is convergent in
$T_{\omega}H_0(S^1,M^n)$ and, thus, converges to the vector-field $\xi^{*} \in
T_{\omega}H_0(S^1,M^n)$ such that $\gamma^{*} = exp(\omega,\xi^{*}) \in
$\newline $ H_0(S^1,M^n)$ and $\|\delta
S_{\varepsilon,\tau}(\gamma^{*})\|=0,$ i.e. the contour $\gamma^{*}$ is
the extremal of the functional $S_{\varepsilon,\tau}$.

Lemma 15 is proved.

{\bf C) Existence of critical points.}

First recall the definition of the Ricci curvature tensor and introduce
its analogue for fields defined by 2-forms $F$.

Let $R^i_{jkl}$ be the Riemann curvature tensor and $\{e_1,\dots,e_n\}$ be
an orthonormal basis in the tangent space $T_xM^n$. Then the Ricci curvature
tensor is defined by the formula
$$
K(u,v) = - \sum^{n}_{\alpha=1} (R(u,e_{\alpha})v,e_{\alpha}).
$$

This definition does not depend on the choice of such a basis because this
tensor can be interpreted  as the trace of the linear transformation
$$
\zeta \rightarrow R(u,\zeta)v.
$$

Let us also consider the tensor
$$
G(u,v,w) = (\nabla_uF)_{ij}v^jw^k.
$$
We can now obtain another tensor
$$
\hat G(u,v,w) = (G(u,v,w) + G(v,u,w))/2
$$
which is symmetric with respect to pair $(u,v)$ and
$$
\hat G(u,u,w) = G(u,u,w)
$$
holds.

Now we introduce the following tensor
$$
H(w) = \sum^{n}_{\alpha=1} \hat G(e_{\alpha},e_{\alpha},w)
$$
in the same manner which we used to introduce the Ricci curvature tensor.
This definition is also independent on the choice of basis $\{e_{\alpha}\}$
because this tensor can be interpreted such as the trace of the linear
transformation
$$
\zeta^i \rightarrow g^{ik} \hat G_{jkl}\zeta^jw^l
$$
where
$$
\hat G_{jkl}u^jv^kw^l = \hat G(u,v,w).
$$

\proclaim{Lemma 16}
Let the function
$$
f(v) = K(v,v)/|v| - H(v)
$$
be positive on the subset of non-zero vectors,
$\gamma$ be a critical point of the functional $S_{\varepsilon,\tau}$ where
$\tau <1$, and $\mu(\gamma)$ be its Morse index.
Then
$$
length(\gamma) \leq \max \{1, \frac{4n((\mu(\gamma)\pi)^2}{\Delta} \}
$$
 where
$$
\Delta = \min_{v \in TM^n, |v|=1}(K(v,v)/|v| - H(v)).
$$
\endproclaim

{\bf Proof of Lemma 16.}

Let take $k > \mu(\gamma)$.

Let $\theta_1(t),\dots,\theta_n(t)$ be a family of
parallel orthonormal vector-fields along $\gamma.$ Since in general
$\theta_i(0) \neq \theta_i(1)$, let us define for every natural number $k$
the family $\beta_{jk}(t)$ by the following formulas
$$
\beta_{jk}(t) = \sin (\pi kt) \theta_j(t)
$$
and
$$
\beta_{jkm}(t) =
\cases
\beta_{jk}(t), &  \frac{m-1}{k} \leq t \leq \frac{m}{k} \\
0, & {\textstyle otherwise}
\endcases.
$$

These vector-fields  lie in
$H_0(S^1,M^n)$ (i.e. in $T_{\gamma}H_0(S^1,M^n)$).

By using elementary computations we obtain that
$$
(\beta_{jkm},\beta_{j'k'm'})_H = \int^{1}_{0}
((\beta_{jkm},\beta_{j'k'm'}) +
(\frac{D}{\partial t}
\beta_{jkm},\frac{D}{\partial t}\beta_{j'k'm'}))dt = 0
$$
if triples $(j,k,m)$ and $(j',k',m')$ do not coincide.

Since
$$
\frac{D}{\partial t} \beta_{jk}(t) = k\pi \cos(\pi kt)\theta_j,
$$
we conclude that
$$
(\frac{D}{\partial t}(\beta_{jk}(t),\frac{D}{\partial t}\beta_{jk}(t)) =
(k\pi)^2\cos^2(\pi kt),
$$
$$
\big(\frac{D}{\partial t}\beta_{jk}(t),\dot \gamma(t)\big) = k\pi \cos(\pi
kt)(\theta_j(t),\dot \gamma(t)),
$$
$$
F_{rs}\beta^r_{jk}\left(\frac{D}{\partial t}\beta_{jk}\right)^s
= \sin{\pi kt} F_{rs} \beta^r_{jk}
\big(\frac{D \theta_{j}}{\partial t}\big)^s +
\pi k F_{rs} \beta^r_{jk} \beta^s_{jk} \equiv 0.
$$

With the help of these formulas let us compute
$$
C_{km}= \sum^{n}_{j=1}
\delta^2S_{\varepsilon,\tau}(\gamma)(\beta_{jkm},\beta_{jkm}).
$$
For the sake of simplicity, let us denote by $L$ the length of the contour
$\gamma$.

By Lemma 6,
$$
C_{km} = 2\varepsilon
\int^{m/k}_{(m-1)/k}(n(k\pi)^2\cos^2(\pi kt) -
K(\dot\gamma,\dot\gamma)) dt +
$$
$$
(1+\tau)L^{\tau-1}
\int^{m/k}_{(m-1)/k}((\tau-1)(k\pi)^2\cos(\pi kt)^2 +
n(k\pi)^2\cos^2(\pi kt) -
$$
$$
K(\dot\gamma,\dot\gamma)\sin^2(\pi kt)) dt +
$$
$$
\int^{m/k}_{(m-1)/k}(\sum^{n}_{j=1}\sin^2(\pi kt)
(\nabla_iF)_{rs}\theta^i_j\theta^r_j\dot\gamma^s)dt =
$$
$$
\frac{2n\varepsilon(k\pi)^2 + (1+\tau)(n-1+\tau)(k\pi)^2L^{\tau-1}}{2k} \ -
$$
$$
(2\varepsilon+(1+\tau)L^{\tau-1})
\int^{m/k}_{(m-1)/k}\sin^2(\pi kt)
K(\dot\gamma,\dot\gamma)dt +
$$
$$
\int^{m/k}_{(m-1)/k}(\sum^n_{j=1}\sin^2(\pi kt)
(\nabla_iF)_{rs}\theta^i_j\theta^r_j\dot\gamma^s)dt.
$$

Since $\varepsilon,\tau < 1$, we conclude that if $L>1$ then
$$
C_{km}
 \leq
\frac{2n(k\pi)^2}{k} - L^{-1}
\int^{m/k}_{(m-1)/k}\sin^2(\pi kt)K(\dot\gamma,\dot\gamma)
dt +
$$
$$
\int^{m/k}_{(m-1)/k}\sin^2(\pi kt)(\sum^n_{j=1}(\nabla_iF)_{rs}
\theta^i_j\theta^r_j\dot\gamma^s),
$$
i.e.
$$
C_{km} \leq 2nk\pi^2 - L^{-1}
\int^{m/k}_{(m-1)/k}\sin^2(\pi kt)K(\dot\gamma,\dot\gamma)
dt + \int^{m/k}_{(m-1)/k}\sin^2(\pi kt)H(\dot\gamma)dt =
$$
$$
2nk\pi^2 -\int^{m/k}_{(m-1)/k}\sin^2(\pi kt)
(K(\dot\gamma,\dot\gamma)/|\dot\gamma| - H(\dot\gamma))dt
$$
$$
\leq 2nk\pi^2 - \frac{L\Delta}{2k}.
$$

Thus, if
$$
L > \max\{1 , \frac{4n(k\pi)^2}{\Delta} \}
$$
then
$$
C_{km} < 0
$$
for every $m=1,\dots,k$.

Therefore for every $m$ there exists $j(m)$ such that
$$
\delta^2 S_{\varepsilon,\tau}(\gamma)(v_m,v_m) <0
$$
where $v_m = \beta_{j(m)km}$.

Now it is easy to see that if

1) $\delta^2 S_{\varepsilon,\tau}(\gamma)(v_i,v_j) = 0$ for $i \neq j$ ;

2) $(v_i,v_j)_H = \delta_{ij}$.

{}From the existence of such a family $\{v_j\}$ it immediately follows that
the Morse index $\gamma$ is not less than $k$ which contradicts
the definition of $\mu(\gamma).$

Lemma 16 is proved.

\proclaim{Lemma 17}
 Let $\{\varepsilon_k,\tau_k\}$ be a sequence of pairs
of positive sufficiently small numbers such that
$$
\lim_{k \rightarrow \infty} (\varepsilon_k,\tau_k) = (0,0)
$$
and $\{\overline{\gamma_k}\}$ be a sequence of points from $H_0(S^1,M^n)$
such that

1) $0 < K_0 \leq E_1(\overline{\gamma_k}) \leq K_1 < +\infty$
where $K_0,K_1$ are some constants which does not depend on $k$;

2) $\delta S_{\varepsilon_k,\tau_k}(\overline{\gamma_k})  =0.$

Then there exists a subsequence $\{\gamma_l\} \subset
\{\overline{\gamma_k}\}$ which converges in $H_0(S^1,M^n)$ to the critical
point of the functional $S=(E_0+\Phi)$:
$$
\lim_{l \rightarrow \infty} \gamma_l = \gamma_{\infty},
$$
$$
\delta (E_0+\Phi)(\gamma_{\infty}) = 0.
$$
\endproclaim

{\bf Proof of Lemma 17.}

By Lemma 4, the contours $\overline{\gamma_k}$ are naturally parameterized
and, since the manifold $M^n$ is compact, we can pick up
a subsequence $\{\gamma_l\}$ such that

1) the lengths of $\gamma_l$ converges to some constant $C>0$:
$$
length(\gamma_l) = C_l \rightarrow C;
$$

2) $\gamma_l(0)$ converges to some point $x_0 \in M^n$;

3) $\dot\gamma_l$ converges to some vector $v \in T_{x_0}M^n$
\newline
as $j \rightarrow \infty$.

By Lemma 5, the motion of the point along every contour
$\gamma_l$ is a periodic solution to the equation
$$
\frac{D}{\partial t}\dot\gamma^j_l = g^{jk}F_{ki}\dot\gamma^i_l/(2\varepsilon_l
+(1+\tau_l)|\dot\gamma_l|^{\tau_l-1}).
$$

It follows from the well-known theorems on continuity
dependence, of solutions of differential equations, on initial data
and coefficients that for every $t$ the sequence of points $\gamma_l(t)$
converges to a point $\gamma_{\infty}(t)$ where the curve
$\gamma_{\infty}$ satisfies the equation
$$
\frac{D}{\partial t}\dot\gamma^j_{\infty} = Cg^{jk}F_{ki}
\dot\gamma^i_{\infty}\eqno{(3.6)}
$$
with the following initial data
$$
\gamma_{\infty}(0)=x_0, \dot\gamma_{\infty}(0)=v.
$$
Since the curves $\gamma_l$ are closed, we may conclude that
the curve $\gamma_{\infty}$ is also closed and, since it satisfies equation
(3.6), this closed curve is a critical point of the functional $S=(E_0+\Phi)$.

Lemma 17 is proved.

\proclaim{Theorem 1}
Let there exists a contour $\gamma \in H_0(S^1,M^n)$
such that
$$
S(\gamma) = \int_{\gamma} (\sqrt{g_{ij}\dot\gamma^i\dot\gamma^j} +
A_i\dot\gamma^i) dt<0.
$$
Let also that
$$
f(v) = K(v,v)/|v| - H(v) > 0\eqno{(3.7)}
$$
for $v\neq 0$.
Then for every cycle $\alpha\in H_l(M^n)\setminus \{0\}$ there exists
a sequence $\{\gamma_{\varepsilon_k,\tau_k}\}$ of critical points of the
functional $S_{\varepsilon_k,\tau_k}$ such that

1) $\varepsilon_k, \tau_k \rightarrow 0$ as $k\rightarrow \infty$;

2) $S(\gamma_{\varepsilon_k,\tau_k}) = c_{\varepsilon,_k,\tau_k}(\alpha)$;

3) $index\ \  \gamma_{\varepsilon_k,\tau_k} \leq (l+1)$;

4) the sequence $\gamma_{\varepsilon_k,\tau_k}$ converges to non-point
extremal of the functional $S$ as $k\rightarrow \infty$.
\endproclaim

{\bf Proof of Theorem 1.}

By Lemmas 9 and 10,
$$
\dim H_{l+1}(H_0(S^1,M^n),\{S_{\varepsilon,\tau} \leq 0\};{\bold R}) \neq 0
$$
for sufficiently small values of $\varepsilon,\tau$.
By Lemma 15, the functionals $S_{\varepsilon,\tau}$
satisfy the Palais-Smale conditions, and, therefore,
there exists  a critical point $\gamma_{\varepsilon,\tau}$ such that
its Morse index is less of equal to
$(l+1)$ and $S_{\varepsilon,\tau}(\gamma_{\varepsilon,\tau}) =
c_{\varepsilon,\tau}(\alpha)$ (the critical values
$c_{\varepsilon,\tau}$ were defined at Lemma 9).

By Lemma 12, there exists a positive constant $K_0$ such that
$$
E_1(\gamma_{\varepsilon,\tau}) \geq K_0,
$$
and, by Lemma 16, there exists the positive constant $K_1$ such that
$$
E_1(\gamma_{\varepsilon,\tau}) = E^2_0(\gamma_{\varepsilon,\tau}) \geq K_1.
$$
Here we ought to take into account that all contours
$\gamma_{\varepsilon,\tau}$ are naturally parameterized (see Lemma 4).

Now it immediately follows from Lemma 17 that there exists a sequence of
$\{\gamma_{\varepsilon_k,\tau_k}\}$ which converges in $H_0(S^1,M^n)$ to a
critical point $\gamma_{\infty}$ of the functional
$S= (E_0+\Phi)$ as $k \rightarrow \infty$ and $E_1(\gamma_{\infty}) >
K_0 >0$ (i.e. the extremal $\gamma_{\infty}$ is not one-point).

Theorem 1 is proved.

Let us notice that the problem of distinguishing between
critical points which we have obtained can be solved under
additional conditions involving the Morse type of
these critical points.

{\bf D) On applications of Theorem 1.}

\proclaim{ Definition}
Let $(M_1,g^{(1)},F^{(1)}),\dots,(M_k,g^{(k)},
F^{(k)})$
be Riemannian manifolds
\newline
endowed with Riemannian metrics $g^{(i)}$ and 2-forms
$F^{(j)}$ respectively. We call the direct product of the Lagrangian systems
$(M_1,g^{(1)},F^{(1)}),\dots,(M_k,g^{(k)},F^{(k)})$ the Lagrangian system of
the same type defined on the direct product of these manifolds
$$
M^{*} = M_1 \times \dots \times M_k
$$
with the product metric and the 2-form $F^{*}$ defined by the following
formula
$$
F^{*} = \sum_{j} p^{*}_j (F^{(j)})
$$
where
$$
p_j: M^{*} \rightarrow M_j
$$
are natural projections on the factors.
\endproclaim

\proclaim{Lemma 18}
 Let $M^n$ be a Riemannian manifold endowed by
an exact non-zero 2-form $F$ such that inequality (3.7) holds.
Then there exists an integer number $k$ such that

1) the product of $k$ examples of the Lagrangian systems $(M^n,g_{ij},F)$
satisfies (3.7) ;

2) on $(M^n)^k$, there exist closed smooth contours
at which the functional
$$
S(\gamma)=\int_{\gamma}(|\dot\gamma|+A^{*}_i\dot\gamma^i)dt
$$
(where $dA^{*}=F^{*}$) takes negative values ,
i.e., this system satisfies the conditions of Theorem 1.
\endproclaim

{\bf Proof of Lemma 18.}
Since the form $F$ is non-zero there exists a loop
$\gamma$, on $M^n$, such that
$$
\int_{\gamma} A_i\dot\gamma^i dt = R<0.
$$

Let $k> (length(\gamma)/R)^2.$ We consider a contour
$\gamma^{*}$ on $(M^n)^k$ which projects onto the contour
$\gamma$ under the projection onto every factor. It is easy to see that
$$
length(\gamma^{*}) = \sqrt{k} \ length(\gamma),
$$
and
$$
\int_{\gamma^{*}} A^{*}_i \dot\gamma^{*i} dt = k
\int_{\gamma}A_i\dot\gamma^i dt.
$$
We then have :
$$
S(\gamma^{*}) = \sqrt{k}\ length(\gamma) - k \int_{\gamma}A_i\dot\gamma^i
dt <0.
$$

Lemma 18 is proved.

{\bf E) On the existence of extremals of the functional S.}

\proclaim{ Theorem 2}
Let $M^n$ be a closed Riemannian manifold endowed by
a Riemannian metric $g_{ij}$ and an exact 2-form $F$ which satisfy the
following condition
$$
\min_{|v|=1}\{f(v) = K(v,v)/|v| - H(v)\}>0.
$$
Then there exists a closed nontrivial extremal of the functional
$$
S(\gamma)=\int_{\gamma}(\sqrt{g_{ij}\dot x^i\dot x^j} + A_i\dot x^i)dt
$$
where $dA=F$.
\endproclaim

{\bf Proof of Theorem 2.}

The real homologies of the pair $(H_0(S^1,M^n),M^n),$ where
$M^n \subset H_0(S^1,M^n),$
are nontrivial (\cite{K1}). Let us pick up some nontrivial cycle
$u \in H_l(H_0(S^1,M^n),$
\newline
$M^n;{\bold R})$.

It is evident that there exists a sequence of pairs of positive numbers
$\{\varepsilon_k,\tau_k\}$ such that

1) $\varepsilon_k, \tau_k \rightarrow 0$ as $k \rightarrow \infty$;

2) the functionals $S_{\varepsilon_k,\tau_k}$ are positive outside of the set
of
one-point contours, or for every of them there exists non-one-point
contour $\gamma_k\in H_0(S^1,M^n)$ such that
$S_{\varepsilon_k,\tau_k}(\gamma_k)=0.$

These functionals $S_{\varepsilon_k,\tau_k}$ satisfy the PS conditions (see
Lemma 15). Thus, if the first case holds for every one of them there exists
non-one-point
extremal with a Morse index which is less or equal to $l$ (this extremal
corresponds, in the usual manner, to the cycle $u$), otherwise for every
one of
them there exists non-one-point extremal with the Morse index equal to $1$
(by Lemma 9).

We can then repeat the arguments of
the proof of Theorem 1 and
derive in this way proof of Theorem 2.

Theorem 2 is proved.

\vskip1.5cm

The second author (I.A.T.) was supported by
the International Science Foundation
(ISF) (grant no. M1E300).

\enddocument